\documentclass[12pt,reqno]{article}
\def\hybrid{\topmargin      0pt
\oddsidemargin 0pt
\headheight 0pt \headsep 0pt
\textwidth 165true mm
\textheight 231true mm
\marginparwidth 0.0in \parskip  0pt plus 1pt \jot = 1.5ex}

\usepackage{amssymb}
\usepackage{amsmath}
\usepackage{amsthm}
\hybrid
\usepackage{amssymb}
\usepackage{amsmath}
\usepackage{amsthm}
\newcommand{\be}[1]{\begin{eqnarray#1}}
\newcommand{\ee}[1]{\end{eqnarray#1}}
\newtheorem{thm}{Theorem}[section]
\newtheorem{propn}[thm]{Proposition}
\newtheorem{lemma}[thm]{Lemma}
\newtheorem{corollary}[thm]{Corollary}

\newtheorem{question}{Question}
\theoremstyle{definition}
\newtheorem{remark}[thm]{Remark}
\newtheorem{example}[thm]{Example}
\newtheorem{definition}[thm]{Definition}

\newcommand{\tp}{\otimes}

\newcommand{\braid}[2]{{#1}$\lower4pt\hbox{$\tp\atop\raise4pt
            \hbox{$\scriptscriptstyle\Ru $}$}${#2}}
\newcommand{\twist}[2]{{#1}${\,\scriptscriptstyle \Ru}\atop\raise9pt\hbox{$\scriptstyle\tp$} ${#2}}
\newcommand{\twistF}[2]{{#1}${\,\scriptscriptstyle \F}\atop\raise9pt\hbox{$\scriptstyle\tp$} ${#2}}

\newcommand{\A}{\mathcal{A}}

\newcommand{\M}{\mathcal{M}}
\newcommand{\La}{\mathcal{L}}
\newcommand{\Ta}{\mathcal{T}}
\newcommand{\Ha}{\mathcal{H}}
\newcommand{\Ru}{\mathcal{R}}

\newcommand{\U}{\mathcal{U}}
\newcommand{\F}{\mathcal{F}}
\newcommand{\C}{\mathbb{C}}
\newcommand{\Z}{\mathbb{Z}}
\newcommand{\K}{\mathbb{K}}
\newcommand{\N}{\mathbb{N}}

\newcommand{\ad}{\mathrm{ad}}

\newcommand{\g}{\mathfrak{g}}

\renewcommand{\l}{\mathfrak{l}}
\newcommand{\h}{\mathfrak h}

\renewcommand{\[}{{[\![}}
\renewcommand{\]}{{]\!]}}

\newcommand{\la}{\lambda}
\newcommand{\ve}{\varepsilon}

\newcommand{\n}{\nonumber          }
\newcommand{\al}{\alpha}
\newcommand{\bt}{\beta}

\newcommand{\si}{\sigma}
\newcommand{\gm}{\gamma}
\newcommand{\id}{\mathrm{id}}
\newcommand{\card}{\mathrm{card}}

\newcommand{\End}{\mathrm{End}}

\newcommand{\Char}{\mathrm{Char}}
\renewcommand{\Re}{\mathrm{Re}}
\renewcommand{\Im}{\mathrm{Im}}

\newcommand{\Tr}{\mathrm{Tr}}
\newcommand{\T}{\mathrm{T}}

\newcommand{\tr}{\triangleright}
\newcommand{\tl}{\triangleleft}

\begin{document}
\title{Method of quantum characters in equivariant quantization}
\author{J. Donin\footnote{This research is partially supported
by the Israel Academy of Sciences grant no. 8007/99-01.}
\hspace{3pt} and A. Mudrov\\[0.1in]
{\small
Department of Mathematics, Bar-Ilan University, 52900 Ramat-Gan, Israel}
}
\date{}
\maketitle
\begin{abstract}
Let $G$ be a reductive Lie group, $\g$ its Lie algebra, and $M$ a $G$-manifold.
Suppose  $\A_h(M)$ is a $\U_h(\g)$-equivariant quantization  of the
function algebra  $\A(M)$ on $M$.
We develop a method of building $\U_h(\g)$-equivariant quantization on
$G$-orbits in $M$ as quotients of $\A_h(M)$. We are concerned with
those quantizations that may be simultaneously  represented as
subalgebras in $\U^*_h(\g)$ and quotients of $\A_h(M)$. It turns out
that they are in one-to-one correspondence with characters of the
algebra $\A_h(M)$. We specialize our approach to the situation
$\g=gl(n,\C)$, $M=\End(\C^n)$, and $\A_h(M)$ the so-called
reflection equation algebra associated with the
representation of $\U_h(\g)$
on $\C^n$. For this particular case, we present in an explicit form
all possible quantizations of this type; they cover symmetric and
bisymmetric orbits. We build a two-parameter deformation family
and obtain, as a limit case, the $\U(\g)$-equivariant quantization of the
Kirillov-Kostant-Souriau bracket on symmetric orbits.
\end{abstract}
\section{Introduction.}
Let $G$ be a reductive Lie group and $\g$ its Lie algebra. Let
$M$ be a right $G$-manifold and $\A_h(M)$  a quantization of the function
algebra $\A(M)$ on $M$. We suppose that the quantization
is $\U_h(\g)$-equivariant, i.e. $\A_h(M)$ is a left $\U_h(\g)$-module
algebra. We consider the problem of "restricting"  $\A_h(M)$ to
the $\U_h(\g)$-equivariant quantization on orbits in $M$.
This means finding an invariant ideal in $\A_h(M)$, a
deformation of the classical ideal specifying an orbit, such that
the quotient algebra will be a flat deformation of
the function algebra on the orbit.
Our principal example is $M=\End(V)$, where $V$ is the underlying
linear space of a finite dimensional representation of $G$.
We consider $\End(V)$ as a right $G$-manifold with respect to
the action by conjugation and  study quantizations on $G$-invariant
sub-manifolds in $\End(V)$.

The problem of equivariant quantization on $\End(V)$, equipped with
the structure of
the adjoint $\U_h(\g)$-module, was considered in \cite{D1}.
It leads to the study of the so-called
reflection equation (RE) algebras,
\cite{KSkl}. Torsion factored out, they are flat deformations of polynomial
functions on the cone $\End^{\Omega}(V)$ of matrices whose tensor square commutes with
the split Casimir ${\Omega}$, the image of the invariant symmetric element from
$\g^{\tp2}$, \cite{DM2}.

It is easy to show that the RE algebra can be restricted to orbits in
$\End^{\Omega}(V)$ in the case $\g=gl(n,\C)$ and $V=\C^n$.
The problem is to describe such a restriction explicitly, i.e.,
to find an appropriate ideal in the RE algebra and prove flatness of the quotient
algebra as a module over $\C[[h]]$.
To this end, we develop a quantization method, confining ourselves to
those quantizations that may be represented as subalgebras in the function
algebra on the quantum group and as quotients of the RE algebra.
Being simultaneously a subalgebra and a quotient algebra of flat deformations
guarantees flatness of the quantization.

In the classical case, every orbit in $M$ is realized as
a subalgebra in  $\A(G)$ and a quotient of $\A(M)$.
The algebra $\A(M)$ is a comodule over the Hopf algebra  $\A(G)$.
A point $a\in M$ defines the map $g\to ag$ from
$G$ onto $O_a$, the  orbit passing through $a$. It
corresponds to a character $\chi^a$  of
the algebra  $\A(M)$, which defines the reversed arrow from
$\A(O_a)$ to  $\A(G)$. Since $O_a\subset M$, the  algebra
$\A(O_a)$ is also a quotient of  $\A(M)$. The idea of our
method is to quantize this picture.
Suppose $\A(M)$ is quantized together with the  $\A(G)$-comodule
structure, so that  $\A_h(M)$ is a comodule over the
dual Hopf algebra $\U^*_h(\g)$.
Suppose there is a character of the algebra $\A_h(M)$ being a deformation of $\chi^a$.
Then, it defines a homomorphism from  $\A_h(M)$ to $\U^*_h(\g)$
and the image of this homomorphism is a deformation
of $\A(O_a)$. We prove that, conversely, every equivariant homomorphism
$\A_h(M)\to \U^*_h(\g)$ is of this form.


We cannot expect to quantize all orbits in such a way. Indeed,
a deformed algebra has a poorer supply of characters than its classical counterpart,
therefore not every orbit, in general, fits our
scheme.
Any character of $\A_h(M)$ corresponds to a point on $M$ where the Poisson bracket
vanishes. An open question is whether every such a point can be quantized
to a character of $\A_h(M)$.

It this paper, we apply our method to the standard quantum
group  $\U_h\bigl(gl(n,\C)\bigr)$ and $M=\End(\C^n)$.
We give full classification of the Poisson brackets on semisimple
orbits of $GL(n,\C)$ that are obtained by restriction from the RE Poisson
structure on $\End(\C^n)$. We describe all $\U_h\bigl(gl(n,\C)\bigr)$-equivariant
quantizations that can be obtained within our approach and present them explicitly
in terms of the RE algebra generators and relations.

In particular, we build quantizations of symmetric and bisymmetric
orbits\footnote{Those are the orbits consisting of matrices with
two and three different eigenvalues, respectively.}
$GL(n,\C)$.
On symmetric orbits, the admissible Poisson brackets form
a one-parameter family and we construct quantizations for all of them.
On bisymmetric orbits, the admissible Poisson brackets are parameterized
by a two dimensional variety. Within our approach, we quantize certain
one-parameter sub-families.
Also, we build quantizations on nilpotent orbits formed by matrices
of zero square.

It is known that a bisymmetric orbit has the structure of a homogeneous
fiber bundle over a symmetric orbit as a base. This fact is important
for the Penrose transformation theory. We show that our quantization
respects that structure and build the quantization of the bundle map.

We extend the quantization on symmetric orbits obtained by the method of characters
to a two-parameter $\U_h\bigl(gl(n,\C)\bigr)$-equivariant family.
It is a restriction of the two-parameter
quantization $\La_{h,t}$ on $\End(\C^n)$, which has the
algebra $\U\bigl((gl(n,\C)\bigr)[t]$
as the limit $h\to 0$. This algebra is a $\U\bigl((gl(n,\C)\bigr)$-equivariant
quantization of the Poisson-Lie bracket on $gl^*(n,\C)\simeq\End(\C^n)$.
Taking the limit $h\to 0$ in the two-parameter deformation
on symmetric orbits, we obtain explicitly the $\U\bigl((gl(n,\C)\bigr)$-equivariant
quantizations of
the Kirillov-Kostant-Souriau (KKS) bracket on symmetric spaces
as quotients of the algebra $\U\bigl(gl(n,\C)\bigr)[t]$.

The setup of the article is as follows.
The next section contains some basic information essential for
our exposition. In particular,
we collect some facts from the quantum group theory in Subsection \ref{ssQGU}
and  recall definitions of modules and comodules over Hopf algebras
in Subsection \ref{ssFRTRE}.
Therein, we introduce the FRT\footnote{The FRT algebra
is the quantization of the Drinfeld-Sklyanin bracket on
$\End(V)$. It was used by Faddeev, Reshetikhin, and Takhtajan
for definition of Hopf algebra  duals to quantum groups.}
and RE algebras associated with the representation of $\U_h(\g)$
on $V$ and recall their basic properties.
In Section \ref{sQAC}, we formulate the method of restricting
the RE algebra to adjoint orbits in $\End(V)$ by means of the RE algebra
characters. We specialize this method to the $GL(n,\C)$-case in
Section \ref{sGL}. We compute the RE Poisson structures on semisimple
(co)adjoint orbits of $GL(n,\C)$ in Subsection \ref{ssREPS} and
present the classification of the RE algebra
characters in Subsection \ref{ssREAC}. On this ground, we build the
quantizations of symmetric and bisymmetric orbits in Subsection \ref{ssQSBO}.
In Subsection \ref{ssQFB}, we show that the
constructed quantization  respects the structure
of homogeneous fiber bundles on bisymmetric orbits.
In Subsection \ref{QKKS}, we construct the two-parameter quantization
on symmetric orbits and give the explicit quantization of the KKS bracket
on them as a limit case.

\section{Preliminaries.}
\label{sAQGS}
\subsection{Quantum group $\U_h(\g)$.}
\label{ssQGU}
Let $\g$ be a reductive Lie algebra over $\C$ and $r\in \wedge^2\g$ a solution to
the modified classical Yang-Baxter equation
\be{}
\[r,r\]=\phi,
\label{phi}
\ee{}
where $\[\cdot,\cdot\]$ stands for the Schouten bracket and $\phi$
is an invariant element from $\wedge^3\g$.
The universal enveloping algebra $\U(\g)$ is a Hopf one, with
the coproduct $\Delta_0$, counit $\ve_0$, and antipode $\gm_0$ defined by
$$
\Delta_0(X) = X\tp 1 +  1\tp X,\quad \ve_0(X)  = 0,  \quad \gm_0(X)=-X,  \quad X\in \g.
$$
These operations are naturally extended over $\U(\g)[[h]]$ as a topological
$\C[[h]]$-module.
The following theorem is implied by the results of Drinfeld, \cite{Dr2}, Etingof and Kazhdan, \cite{EK},
\begin{thm}
There exists an element $\F_h  \in \U^{\tp 2}(\g)[[h]]$,
$$
\F_h=1\tp 1 + \frac{h}{2} r + o(h),
$$
such that $\U(\g)[[h]]$ becomes a quasitriangular Hopf algebra  $\U_h(\g)$ with
the coproduct $\Delta$, counit $\ve$, and antipode $\gm$:
\be{}
\Delta(x) &=& \F_h^{-1} \Delta_0(x) \F_h
,\n\\
\ve(x) &=& \ve_0(x)
,
\label{drtw}
\\
\gm(x) &=& u^{-1}\gm_0(x)u
,\quad
 x \in \U_h(\g)
.\n
\ee{}
The element $u$ is equal to $ u = \gm_0(\F_{1})\F_{2}$, where $\F_1\tp \F_2=\F_h$ (summation
implicit) and
the universal R-matrix is given by
\be{}
\label{urm}
\Ru_h = (\F_h^{-1})_{21}e^{h\omega} \F_h = 1\tp 1 + h\Bigl(r+\omega\Bigr) + o(h),
\ee{}
where $\omega \in \g^{\tp 2}$ is a symmetric invariant element
such that the sum $r+\omega$ satisfies the classical Yang-Baxter
equation, \cite{Dr1}.
\end{thm}
The algebra $\U_h(\g)$ is a quantization of $\U(\g)$ in the sense that
$\U_h(\g)/h\:\U_h(\g)=\U(\g)$ as a Hopf algebra.
The coassociativity of the coproduct $\Delta$ implies the cocycle equation
on $\F_h$:
\be{}
(\Delta\tp id)(\F_h)(\F_h\tp 1)=\Phi^{-1}_h(id \tp \Delta)(\F_h)(1\tp\F_h).
\label{cocycle}
\ee{}
Here, $\Phi_h$ is an invariant element from $\U^{\tp 3}_h(\g)$;
it is called {\em coassociator} and satisfies the pentagon identity
$$
(\id^{\tp 2} \tp \Delta)(\Phi_h)(\Delta\tp \id^{\tp 2} )(\Phi_h)=
(1\tp\Phi_h)(\id \tp \Delta \tp \id)(\Phi_h)(\Phi_h\tp 1).
$$
\begin{example}[Standard quantum groups]
Let $\Pi$ be the root system of $\g$ and $\Pi^\pm$
the subsets of positive and negative roots.
Let $e_{\pm\al}$, $\al \in \Pi^+$, be root vectors normalized
to $(e_{\al},e_{-\al})=1$ with respect to the Killing form.
The simplest example of the classical r-matrix is
\be{}
r=\sum_{\al\in \Pi^+} e_\al \wedge e_{-\al}.
\label{str-m}
\ee{}
It is called {\em standard} r-matrix and the corresponding
quantum group {\em standard} or Drinfeld-Jimbo quantization of $\U(\g)$.
Other possible r-matices for simple Lie algebras are listed
in \cite{BD}. They were explicitly quantized in \cite{ESS}.
\end{example}

By $\U^*_h(\g)$ we mean the FRT dual to $\U_h(\g)$, \cite{FRT}.
This is a quantized polynomial algebra on the group $G$;
as a linear space, $\U^*_h(\g)$ consists of
$\End^*(V)$ while $V$ runs over finite dimensional completely reducible
representations of  $\U_h(\g)$.

\subsection{FRT and RE algebras.}
\label{ssFRTRE}
In this section, we collect some basic facts about the
FRT and RE algebras and their symmetries.
All $\U_h(\g)$-modules are considered to be free over
$\C[[h]]$. Tensor products are assumed to be completed in the $h$-adic topology.

Recall that an associative algebra $\A$ over $\C[[h]]$ is called a
left (right)
$\U_h(\g)$-{\em module algebra} if it is a left (right) module with
respect to the action $\tr$ ($\tl$) and its multiplication
is consistent with the module structure:
\be{}
x\tr(a b) = \bigl(x_{(1)}\tr a\bigr) \bigl(x_{(2)}\tr b\bigr),
& &
(a b)\tl x= \bigl(a\tl x_{(1)}\bigr) \bigl(b\tl x_{(2)}\bigr),
\label{lractions}
\\
1\tr a = a,
& &
a\tl 1= a,
\label{lrunit}
\\
x\tr 1_\A = \ve(x)1_\A,
& &
1_\A\tl x= \ve(x)1_\A
\label{unitlr}
\end{eqnarray}
for any $x\in \U_h(\g)$ and  $a,b\in \A$.
We adopt the standard brief Sweedler notation for the coproduct
$\Delta(x)=x_{(1)}\tp x_{(2)}$, where $x$ is an element from a Hopf
algebra $\Ha$.

If $\A$ is a left and right
module  simultaneously and the two actions commute
with each other,
\be{}
\label{l-r}
x_1\tr (a \tl x_2)= (x_1\tr a) \tl x_2 ,
\quad
x_1,x_2\in \U_h(\g), \> a\in \A,
\ee{}
then it is called {\em bimodule}.
$\A$ is a $\U_h(\g)$-{\em bimodule algebra}
if its bimodule and algebra structures are
consistent with the coproduct in $\U_h(\g)$
in the sense of  (\ref{lractions}--\ref{unitlr}).

A right $\U^*_h(\g)$-{\em comodule algebra} is an associative algebra $\A$
endowed with a homomorphism $\delta\colon \A\to \A\tp \U^*_h(\g)$ obeying
the coassociativity constraint
\be{}
(\id\tp \Delta)\circ\delta &=& (\delta\tp \id)\circ\delta
\label{coass}
\ee{}
and the conditions
\be{}
\delta(1_\A)=1_\A\tp 1,\quad (\id\tp \ve)\circ\delta  =  \id,
\label{coun}
\ee{}
where the identity map on the right-hand side of the
second equation assumes the isomorphism
$\A\tp \C[[h]]\simeq \A$.
As for the coproduct $\Delta$, we use symbolic notation
$\delta(x)=x_{[1]}\tp x_{(2)}$, marking the tensor component
belonging to $\A$ with the square brackets. The subscript
of the $\U^*_h(\g)$-component is concluded in parentheses.
Every  right $\U^*_h(\g)$-comodule $\A$ is a left $\U_h(\g)$-module,
the action being defined through the pairing $\langle \:.\: ,\!\:.\: \rangle$
between $\U_h(\g)$ and $\U^*_h(\g)$:
\be{}
\label{coaction}
x\tr a &=& a_{[1]}\langle x ,a_{(2)} \rangle,
\quad x\in \U_h(\g),\> a\in  \A.
\ee{}
If $\A$ is finite dimensional, the converse is also true.
Similarly to right $\U^*_h(\g)$-comodule algebras, one can
consider left ones. They are also right $\U_h(\g)$-module algebras.

A completely reducible finite dimensional representation $\rho$ of the universal enveloping
algebra $\U(\g)$ on a linear space $V$ is naturally extended to that of
$\U_h(\g)$ on $V[[h]]$. Denote by $\M$ the matrix space $\End(V)$
and fix a basis $e_i$ in $V$. Let $e^i_j \in \M$
be the matrix units acting on $V$ from the right by the rule
$ e_l e^i_j = \delta_l^i e_j$, where $\delta^i_l$ is the Kronecker symbol.
In terms of the basis $\{e^i_j\}$, the multiplication is expressed by the formula
$e^i_j e^l_k = \delta^l_j e^i_k$.

As an associative algebra, $\M$ is a bimodule over itself with respect
to the right and left regular representation. The homomorphism $\rho$ equips
$\M[[h]]$ with the structure of a  $\U_h(\g)$-bimodule. By duality, the space
$\M^*[[h]]$ is a  $\U_h(\g)$-bimodule as well.

Let $R$ denote the image of the universal R-matrix of $\U_h(\g)$
under the representation $\rho$. We shall also use the matrix $S$, which
differs from $R$ by the permutation $P$ on $V\tp V$,
\be{}
\label{S=PR}
S=PR,\quad P = \sum_{i,j}e^i_j\tp e^j_i.
\ee{}
\begin{example}[FRT algebra]
Let $\{T^i_k\}\subset \M^*$ be the dual basis to $\{e^k_i\}$.
The associative algebra $\Ta_R$ is generated by the matrix
coefficients $\{T^i_k\} \subset \M^*$ subject to the FRT relations
\be{}
\label{frt0}
\sum_{\al,\bt} S_{ij}^{\al\bt} T^m_\al T^n_\bt  =
\sum_{\al,\bt}  T^\al_i T^\bt_j S_{\al\bt}^{mn},
\ee{}
or, in the standard compact form,
\be{}
ST_1 T_2  =  T_1T_2 S,
\label{frt}
\ee{}
where $T=\sum_{i,j} T^i_j e^j_i$.
The matrix elements $T^i_j$  of the representation $\rho$
may be thought of as linear functions on  $\U_h(\g)$; they define
an algebra homomorphism
\be{}
\label{TtoH}
\Ta_R\to\U^*_h(\g).
\ee{}
\end{example}
\begin{propn}
\label{FRTA}
Let $\rho$ be a finite dimensional  completely reducible representation of $\U_h(\g)$ on
a module $V[[h]]$ and $\Ta_R$ the FRT algebra
associated with $\rho$. Then,
$\Ta_R$
is a $\U_h(\g)$-bimodule algebra, with the left and right actions
extended from $\M^*$:
\be{}
x\tr T = T\rho(x),\quad   T\tl x  = \rho(x) T, \quad x \in \U_h(\g).
\label{lra}
\ee{}
It is a bialgebra, with the coproduct and counit being defined as
\be{}
\Delta(T^i_j)=\sum^n_{l=1}T^l_j\tp T_l^i,
\quad
\ve(T^i_j) =\delta^i_j.
\ee{}
Composition of the coproduct with the algebra homomorphism (\ref{TtoH})
applied to the (left) right tensor factor makes
$\Ta_R$ a (left) right $\U^*_h(\g)$-comodule
algebra.
\end{propn}
\begin{proof}
Actions (\ref{lra}) are extended to the actions on the tensor
algebra $\T(\M^*)[[h]]$, and the ideal generated by
(\ref{frt}) turns out to be invariant. Concerning the bialgebra structure,
the reader is referred to \cite{FRT}. The structure of a comodule is
inherited
from the bialgebra one, so it is obviously coassociative.
The coaction is an algebra homomorphism, being a composition of
two homomorphisms.
\end{proof}
\noindent
Remark that the FRT relations (\ref{frt}) arose within the quantum inverse
scattering method and were used for systematic definition of the
quantum group duals in \cite{FRT}.
\begin{example}[RE algebra]
Another algebra of our interest, $\La_R$, is
defined as the quotient of the tensor algebra $\T\bigl(\M^*)[[h]]$ by the RE relations
\be{}
\label{re0}
\sum_{\al,\bt,\mu,\nu}
S_{ij}^{\al\bt} L^\mu_\bt S^{m\nu}_{\al\mu} L_\nu^n
=
\sum_{\al,\bt,\mu,\nu}
L^\al_j S^{\bt\mu}_{i\al} L^\nu_\mu S^{mn}_{\bt\nu},
\ee{}
where $\{L^i_j\}$ is the basis in $\M^*$ that is dual to
$\{e^j_i\}$. In the compact form, (\ref{re0}) reads
\be{}
SL_2 S L_2  =  L_2 S L_2 S,
\label{re}
\ee{}
where $L$ is the matrix $\sum_{i,j} L^i_j e^j_i$.
\end{example}
\begin{propn}
\label{REA}
Let $\rho$ be a finite dimensional  completely reducible representation of $\U_h(\g)$ on
a module $V[[h]]$ and $T^k_l\in \U^*_h(\g)$
its  matrix coefficients.
Let $L^i_j$ be the generators of the algebra $\La_R$ associated with $\rho$.
Then, $\La_R$ is a left $\U_h(\g)$-module algebra with the action
extended from the coadjoint representation in $\M^*[[h]]$:
\be{}
x\tr L = \rho\bigl(\gamma(x_{(1)})\bigr)L\rho(x_{(2)})), \quad x \in \U_h(\g).
\label{ada}
\ee{}
It is a right $\U^*_h(\g)$-comodule algebra with respect to the coaction
\be{}
\delta(L^i_j) = \sum_{l,k} L^l_k\tp \gamma(T^k_j) T^i_l,
\label{ada-re}
\ee{}
\end{propn}
\begin{proof}
Action (\ref{ada}) is naturally extended to  $\T(\M^*)[[h]]$ and it preserves
relations (\ref{re}). The coassociativity of (\ref{ada-re}) is obvious.
To prove that $\delta$ is an algebra homomorphism, one needs to employ
commutation relations $(\ref{frt})$ and $(\ref{re})$. For details,
the reader is referred to \cite{KS}.
\end{proof}
\noindent
A spectral dependent version of the RE appeared first in
\cite{Cher}. In the form of (\ref{re}), it may be found
in articles \cite{Skl,AFS} devoted to integrable models. The algebra
$\La_R$ was studied in \cite{KSkl,KS}. Its relation to the braid
group of a solid handlebody was
pointed out in \cite{K}.

Remark that the algebras $\Ta_R$ and $\La_R$ may be defined
for any quasitriangular Hopf algebra and its finite dimensional
representation. Propositions \ref{FRTA} and \ref{REA} will be
also valid.

It was proven in  \cite{DM2} that the FRT and RE algebras
are twist-equivalent. For a detailed exposition of the
twist theory, the reader is referred e.g. to \cite{Mj}.
Here we recall that the twist of a Hopf
algebra $\Ha$ with a cocycle $\F$ is a Hopf
algebra $\tilde \Ha$ with the same multiplication as in
$\Ha$ and the coproduct
\be{}
\label{tcopr}
\tilde \Delta(x)=\F^{-1}\Delta(x)\F, \quad x\in \Ha.
\ee{}
To ensure coassociativity of $\tilde \Delta$, the element $\F$ must obey
the constraint
\be{}
(\Delta \tp \id)(\F)\F_{12} = (\id\tp \Delta)(\F)\F_{23}.
\label{cc}
\ee{}
If $\A$ is a left $\Ha$-module algebra with the multiplication $m$,
the new associative multiplication
\be{}
\label{twmod}
\tilde m (a\tp b) = m (\F_1\tr a\tp \F_2\tr b), \quad a, b \in \A,
\ee{}
can be introduced on $\A$. This algebra, $\tilde \A$,
is an $\tilde \Ha$-module algebra.

For example, if $\Ha$ is quasitriangular with the universal
R-matrix $\Ru$, the coopposite
Hopf algebra  $\Ha^{op}$ is a twist with $\F=\Ru^{-1}$. Another
example is the twisted tensor square \twist{$\Ha$}{$\Ha$}
of a quasitriangular Hopf algebra. This is a twist of the
ordinary tensor square $\Ha\tp\Ha$ by the cocycle
$\F=\Ru_{23}\in (\Ha\tp\Ha)\tp(\Ha\tp\Ha)$.
\begin{thm}[\cite{DM2}]
\label{FRT-RE}
Let $\Ta_R$ and $\La_R$ be respectively the FRT and RE algebras
associated with a finite dimensional representation of $\Ha$.
Consider $\Ta_R$ as an $\Ha^{op}\tp\Ha$-module algebra.
Then, there exists a twist from
$\Ha^{op}\tp\Ha$ to \twist{$\Ha$}{$\Ha$} such
that the induced transformation (\ref{twmod}) converts $\Ta_R$ to $\La_R$.
\end{thm}
This twist is performed in two steps. The first one transforms
the factor $\Ha^{op}$ to $\Ha$ in $\Ha^{op}\tp\Ha$.
It is carried out via the
cocycle $\F'=\Ru_{13}\in (\Ha^{op}\tp\Ha)\tp(\Ha^{op}\tp\Ha)$.
The second twist from $\Ha\tp\Ha$ to \twist{$\Ha$}{$\Ha$}
 is via the cocycle
$\F''=\Ru_{23}\in (\Ha\tp\Ha)\tp(\Ha\tp\Ha)$.
The composite transformation with the cocycle
$\F=\F'\F''$ converts multiplication in $\Ta_R$ to
that in $\La_R$ according to formula (\ref{twmod}).

 It follows from Theorem \ref{FRT-RE} that the algebras $\Ta_R$ and $\La_R$ are isomorphic
as $\C[[h]]$-modules in the case $\Ha=\U_h(\g)$.
Another consequence is that $\La_R$ is a module not only
over $\Ha$ but over \twist{$\Ha$}{$\Ha$} as well.
The action of $\Ha$ on $\La_R$ is induced by the Hopf algebra
embedding $\Ha\to \mbox{\twist{$\Ha$}{$\Ha$}}$ via the coproduct.

\section{Quantum characters and quantization on orbits.}
\label{sQAC}
\subsection{Invariant monoid $\M^\Omega$.}
Let  $V$ be a complex linear space and $\rho$ a homomorphism of  $\U(\g)$ into the matrix
algebra $\M=\End(V)$. An element $\xi\in \g$ generates the left and right
invariant vector fields on $\M$
$$
(\xi^l \tr f)(A)=df\bigl(A\rho(\xi)\bigr),
\quad
(f\tl\xi^r)(A) =df\bigl(\rho(\xi)A\bigr),\quad A\in \M, \> f\in \A(\M),
$$
defining the left and right actions of the algebra $\U(\g)$ on
functions on $\M$. Given an element $\psi\in \U(\g)$, by $\psi^r$ and $\psi^l$
we denote, correspondingly, its extensions to the right- and left-invariant
differential operators on  $\M$.
The left adjoint action of $\U(\g)$ on $\M$ is generated by
the vector fields
\be{}
\xi^{\ad}=\xi^l - \xi^r,\quad \xi\in \g .
\ee{}
Let $\Omega\in \M^{\tp2}$ be the image of the invariant symmetric tensor
$\omega\in \g^{\tp2}$ participating in construction of  $\U_h(\g)$, cf.
formula (\ref{urm}).
Introduce  the cone of
matrices\footnote{It coincides with $\M=\End(V)$
for $\g=sl(n,\C)$ and  $V=\C^n$.}
\be{}
\label{cone}
\M^{\Omega}=\{ A\in \M|\; [\Omega , A\tp A]=0\}.
\ee{}
Evidently, $\M^{\Omega}$
is an algebraic variety, it is closed under the matrix multiplication and invariant
with respect to the two-sided action of the group
$G$.
 There are two remarkable
Poisson structures on $\M^\Omega$; they are given by the Drinfeld-Sklyanin bracket
\be{}
r^{l,l}-r^{r,r}
\label{dsbr}
\ee{}
and the RE bracket
\be{}
r^{\ad,\ad}+(\omega^{r,l}-\omega^{l,r}) .
\label{rebr}
\ee{}

\begin{thm}
\label{DSDM1}
\begin{enumerate}
\item
The quotient of  $\Ta_R$ by torsion
is a $\U_h(\g)^{op}\tp\U_h(\g)$-equivariant quantization of
Poisson bracket (\ref{dsbr})
on the cone $\M^{\Omega}$.
\item
The quotient of $\La_R$ by torsion is a $\U_h(\g)$-equivariant quantization of  Poisson bracket  (\ref{rebr})
on the cone $\M^{\Omega}$.
\end{enumerate}
\end{thm}
\begin{proof}
For the proof of the first statement, the reader is referred to
\cite{DS}. The second statement is deduced from the first one
using the twist from  Theorem \ref{FRT-RE}, see \cite{DM2}.
\end{proof}

\subsection{General formulation of the method.}
In the previous section, we considered two examples of equivariant quantization on
the space $\M^{\Omega}$. Depending on a particular choice of symmetry,
they were quotients of the algebras $\Ta_R$ and $\La_R$ by the ideal of $h$-torsion
elements.
Further we study $\U_h(\g)$-invariant ideals in $\La_R$  that are  deformations of
the classical ideals in the function algebra $\A(\M^{\Omega})$
specifying the orbits. The problem is
to ensure flatness of the quotient algebras.
In this section, we formulate a method realizing the quantized
orbits simultaneously as quotients and subalgebras of flat deformations,
and that guarantees flatness of the quantization.
Briefly, we construct quantizations that are
quotients of $\La_R$ and subalgebras in $\A_h(G)=\U^*_h(\g)$.

Our quantization method uses
analogs of points, which are one-dimensional representations or characters
of quantized algebras. Let $M$ be a manifold with a right action
of the group $G$. Let $a\in M$ be a point and $\chi^a$ the
corresponding character of the function algebra
$\A(M)$.
On the diagram
$$
\begin{array}{rcc}
M \times G&\longrightarrow &M\\[8pt]
      \uparrow \hspace{14pt} &&\uparrow\\[8pt]
\{a\}\times G&\longrightarrow &O_a
\end{array}
\quad
\begin{array}{rcc}
\A(M)\tp\A(G)  &\longleftarrow &\A(M)\\[8pt]
\chi^a\tp \id\downarrow  \quad\quad\>\>\; &&\downarrow\\[8pt]
   \C \tp\A(G) &\longleftarrow &\A(O_a)
\end{array},
$$
the left square  displays embedding of the orbit $O_a$
passing through $a$ into the $G$-space $M$.
The induced morphisms of the function algebras are  depicted
on the right square. Our goal is to quantize this picture.

\begin{propn}
\label{coequivar}
Let $\Ha$ be a Hopf algebra over $\C[[h]]$ and $\A$ a comodule algebra over its
Hopf dual $\Ha^*$. Any character $\chi$ of $\A$
defines a homomorphism $\varphi_\chi\colon \A\to \Ha^*$
fulfilling the equivariance condition
\be{}
(\varphi_\chi\tp \id)\circ \delta = \Delta \circ \varphi_\chi.
\label{equivar}
\ee{}
Conversely, any homomorphism $\varphi\colon\A\to \Ha^*$ obeying (\ref{equivar})
has the form $\varphi_\chi$, where $\chi$ is a character of $\A$.
\end{propn}
\begin{proof}
Let $\chi$ be a character of $\A$.
We set $\varphi_\chi=(\chi\tp \id) \circ \delta$ and make
use of the isomorphism $\C[[h]]\tp\Ha^*  \simeq \Ha^*$. Due to coassociativity
(\ref{coass}) of the coaction $\delta$, condition (\ref{equivar}) holds.
Conversely, if $\varphi$ satisfies (\ref{equivar}), then we apply
the counit of the Hopf algebra $\Ha^*$ to the first
tensor factor of (\ref{equivar}) and obtain
$\varphi_\chi$ with $\chi=\ve\circ \varphi$.
\end{proof}

Given an associative algebra $\A$ over $\C[[h]]$ let $\Char(\A)$ denote its set of
characters, i.e. homomorphisms $\A\to \C[[h]]$.
Any element $\chi\in\Char(\A)$ defines an algebra
$\A_\chi$ closing up the commutative diagram (the right-most arrow is onto)
\be{}
\begin{array}{rcc}
\A\tp \Ha^*  &\stackrel{\delta}{\longleftarrow} &\A\\[8pt]
\chi\tp \id\downarrow  \hspace{19pt}&&\downarrow\\[8pt]
\C[[h]] \tp \Ha^*&\longleftarrow &\A_\chi
\end{array}.
\ee{}
We consider $\Ha^*$ as the left coregular module over $\Ha$.
Then, because of (\ref{equivar}),
$\varphi_\chi$ is an $\Ha$-equivariant homomorphism of algebras,
and its image $\A_\chi$ is an  $\Ha$-module algebra.
\begin{definition}
Let $\A$ be a right $\Ha^*$-comodule algebra. Two characters
$\chi_1,\chi_2\in \Char(\A)$ are called gauge-equivalent, $\chi_1\sim \chi_2$,
if there exists an element $\eta\in \Char(\Ha^*)$ such that
\be{}
\label{chareq}
\chi_2 = (\chi_1\tp \eta)\circ \delta .
\ee{}
\end{definition}
This is an equivalence relation. Indeed,
 $\chi_1\sim \chi_2$ and $\chi_2\sim \chi_3$ implies
$\chi_1\sim \chi_3$, due to coassociativity of the coaction.
Also, $\chi_1\sim \chi_2\Rightarrow \chi_2\sim \chi_1$, since
(\ref{chareq}) implies
$\chi_1 = (\chi_2\tp \eta\circ\gm )\circ \delta$. Obviously
$\chi_1\sim \chi_1$ via the counit $\ve\in \Char(\Ha^*)$ of
the Hopf algebra $\Ha^*$.
In the classical situation $\A=\A(M)$ and $\Ha^*=\A(G)$, two characters are gauge-equivalent
if and only if the corresponding points lie on the same orbit
of $G$.
\begin{propn}
Let $\A$ be a right $\Ha^*$-comodule algebra and
$\chi_1\sim \chi_2\in \Char(\A)$. There exist Hopf algebra
automorphisms
$f_{\Ha}\colon \Ha\to\Ha$,
$f_{\Ha^*}\colon \Ha^*\to\Ha^*$
and an algebra automorphism
$f_{\A}\colon \A\to\A$ such that
the diagram
$$
\begin{array}{rcl}
\A&\stackrel{\varphi_{\chi_1}}{\longrightarrow}&\Ha^*\\
f_{\A}\downarrow& &\>\downarrow f_{\Ha^*}\\
\A&\stackrel{\varphi_{\chi_2}}{\longrightarrow}&\Ha^*
\end{array}
$$
is commutative and $\Ha$-equivariant with respect to the left
$\Ha$-actions
$$
\begin{array}{rcl}
\Ha\tp\A&\longrightarrow&\A\\
f_{\Ha}\tp f_{\A}\downarrow\hspace{14pt} & &\>\downarrow f_{\A}\\
\Ha\tp\A&\longrightarrow&\A
\end{array}\quad
\begin{array}{rcl}
\Ha\tp\Ha^*&\longrightarrow&\Ha^*\\
f_{\Ha}\tp f_{\Ha^*}\downarrow\hspace{19pt}& &\>\downarrow f_{\Ha^*}\\
\Ha\tp\Ha^*&\longrightarrow&\Ha^*
\end{array}
$$
\end{propn}
\begin{proof}
Let $\eta$ be the character of $\Ha^*$ realizing the
equivalence $\chi_1\sim \chi_2$ in (\ref{chareq}).
It may be thought\footnote{This is obviously true for finite dimensional Hopf algebras.
In the quantum group case, $\Ha$ is complete in the $h$-adic topology. Then $\Ha^*$ consists of
continuous linear functionals on $\Ha$. It is
equipped with the weak topology, in which continuous linear functionals on $\Ha^*$ form $\Ha$.} of as a group-like element
from $\Ha$, i.e. the one whose coproduct is equal to
$\Delta(\eta)=\eta\tp\eta$.
The similarity transformation with such an element is
an automorphism of the  Hopf algebra $\Ha$.
We set
\be{}
f_{\Ha^*}(x)&=&\eta\bigl(x_{(1)}\bigr)x_{(2)}\eta\Bigl(\gm\bigl(x_{(3)}\bigr)\Bigr)
,\quad x\in \Ha^*,
\label{f_H^*}\\
f_{\Ha}(y)&=&\gm(\eta)\:y\:\eta
,\quad y\in \Ha,\\
f_{\A}(a)&=& a_{[1]}\eta\Bigl(\gm \bigl(a_{(2)}\bigr)\Bigr)
,\quad a\in \A.
\ee{}
A straightforward verification using coassociativity of
$\delta$ and $\Delta$ shows that these maps possess the
required properties.
\end{proof}

Specifically for the deformation quantization situation,
we suppose that $\A_h(M)$ is a quantization of $\A(M)$ and
it is a comodule algebra for
$\A_h(G)\simeq\U^*_h(\g)$. The coaction $\delta\colon \A_h(M)\to\A_h(M)\tp\U^*_h(\g)$
is assumed to be a deformation of the classical map
$\A(M)\to\A(M)\tp\A(G)$. This implies that $\A_h(M)$ is
an equivariant quantization of $\A(M)$ because every
$\U^*_h(\g)$-comodule is a $\U_h(\g)$-module via action (\ref{coaction}).
In connection with Proposition \ref{coequivar}, there arises the problem of
describing the set $\Char\bigl(\A_h(M)\bigr)$.

The following statement is
elementary.
\begin{propn}
\label{clchar}
Let $M$ be a Poisson manifold and $\A_h(M)$ a
quantization of the function algebra $\A(M)$.
Let $a\in M$ and suppose $\chi_h$ is a character of the algebra  $\A_h(M)$
such that $\chi_h(f) = f(a)\! \mod h$, $f\in \A(M)$.
Then, the  Poisson bracket vanishes
at the point $a$.
\end{propn}
\begin{proof}
By definition, $\chi_h\bigl(m_h(f,g)\bigr)=\chi_h(f)\chi_h(g)$,
for $f,g,\in \A(M)$.
Expanding this equality in the deformation parameter and
collecting the terms before $h$, we come to the condition
\be{}
\label{ch0}
m_1(f,g)(a) + \chi_1 (fg) =
f(a) \chi_1 (g) +  \chi_1 (f)g(a).
\ee{}
Here $m_1$ and $\chi_1$ are the infinitesimal terms of the deformed
multiplication $m_h$ and character $\chi_h$.
Skew-symmetrization of (\ref{ch0}) proves the statement.
\end{proof}
Given a Poisson manifold $M$,
we denote by $\Char_0(M)$ the subset of points where the Poisson
bracket vanishes.
\begin{question}
\label{Problem}
Let $M$ be a Poisson manifold and $\A_h(M)$ a quantization of
its function algebra. Given a point $a\in \Char_0(M)$,
does there exist a character $\chi_h \in \Char\bigl(\A_h(M)\bigr)$
such that $\chi_h(f)=f(a)\!\mod h$ for all $f\in \A(M)$?
\end{question}
\subsection{Application to the quantum matrix algebra $\La_R$.}
We specialize the construction suggested in the previous
subsection, to the situation when $\Ha$ is the quantum group $\U_h(\g)$,
$\Ha^*$ the quantized function algebra $\A_h(G)\simeq\U^*_h(\g)$
on the group $G$,
and $M$ the matrix cone $\M^\Omega$ relative to
a given representation of $\U_h(\g)$.
The equivariant quantization $\A_h(M)$ is the RE algebra $\La_R$.
It is a $\U^*_h(\g)$-comodule algebra, and we may apply
Proposition \ref{coequivar} in order to
obtain quatnization of orbits  in $\M^\Omega$
as quotients of $\La_R$.
Elements of $\Char(\La_R)$ are defined by matrices
$A\in \M[[h]]$ solving the numerical reflection equation
\be{}
\label{nre}
S A_2 S A_2 =  A_2 S  A_2 S.
\ee{}
\begin{definition}
\label{conj}
We say that a matrix $A_h\in \M[[h]]$ belongs to the orbit
${O_A}$ if $A_h=A\tl u $ for some invertible
element $u\in \U(\g)[[h]]$.
\end{definition}
\begin{thm}
\label{thmquan}
Let $A_h\in \M[[h]]$ be a solution of (\ref{nre}) and
$\chi=\chi^{A_h}$ the corresponding character of the algebra $\La_R$.
Suppose the matrix $A_h$ belongs to $O_A\subset \M^{\Omega}$. Then, the
algebra  $\A_\chi$ closing up the commutative diagram (the right-most arrow is onto)
\be{}
\label{qd}
\begin{array}{rcc}
\La_R\tp \U^*_h(\g)  &\longleftarrow &\La_R\\[8pt]
\id\tp \chi^{A_h}\downarrow \hspace{30pt} \:&&\downarrow\\[8pt]
\C[[h]] \tp \U^*_h(\g)&\longleftarrow &\A_\chi
\end{array}
\ee{}
is a $\U_h(\g)$-equivariant quantization of the
polynomial algebra on $O_A$.
Its embedding into the left coregular $\U_h(\g)$-module $\U^*_h(\g)$ is equivariant.
\end{thm}
\begin{proof}
The algebra $\A_\chi$ is simultaneously defined by the commutative diagram
(\ref{qd}) as a quotient and a subalgebra of flat $\C[[h]]$-modules, therefore it is
flat. $\A_\chi$ is invariant under the left action
$x\tr T = T\rho(x)$, $x\in\U_h(\g)$, and we should
check that it is isomorphic to  $\A(O_A)[[h]]$ as a $\U(\g)[[h]]$-module.
Let $u$ be an element from $\U_h(\g)$, such that
$A_h=A \tl u$.
The coboundary twist of $\U_h(\g)$ with the element
$\Delta(u)(u^{-1}\tp u^{-1})$
converts $\A_\chi$ into  $\tilde\A_\chi$ for
which $A$ is a character. The matrix $A$ defines an equivariant embedding
$\tilde\A_\chi$ into the twisted algebra
$\tilde\U^*_h(\g)$, which is
generated by the matrix elements of $\tilde T = \rho(u) T \rho(u^{-1})$.
Its image is a subalgebra $\tilde\A_\chi$ in $\tilde\U^*_h(\g)$.
Since
$\U_h(\g)$ itself is a twist of $\U(\g)[[h]]$, the algebras
$\tilde\A_\chi$ and $\A_\chi$ are isomorphic
as $\U(\g)[[h]]$-modules.
Obviously, the subalgebra $\tilde\A_\chi$ coincides with $\A(O_A)$
modulo $h\A(O_A)$.
\end{proof}
The gauge-equivalence between characters of $\La_R$
are realized by means of elements from $\Char\bigl(\U^*_h(\g)\bigr)$, which are
described by the following proposition.
\begin{propn}
\label{char-standard}
For the Drinfeld-Jimbo quantum group $\U_h(\g)$,
the set $\Char\bigl(\U^*_h(\g)\bigr)$ consists of
the elements $e^\eta \in \U_h(\g)$, where
$\eta$ belongs to the Cartan subalgebra $\h\subset \g$.
\end{propn}
\begin{proof}
The standard r-matrix (\ref{str-m}) is of zero weight, so
the subset $\Char_0(G)$ coincides with the maximal torus corresponding
to the Cartan subalgebra.
As an associative algebra,  $\U^*_h(\g)$ is isomorphic to
$\U(\g^*)[[h]]$, \cite{Dr1}. Its characters are parameterized by the dual space to
$\g^*/[\g^*,\g^*]=\h^*$.
On the other hand, the elements $e^\eta$, $\eta \in \h$, are group-like
because all $\eta$ are primitive with respect to the coproduct in
$\U_h(\g)$, \cite{Dr1}.
\end{proof}

In the next section we specialize our consideration to the
Drinfeld-Jimbo quantum group $\U_h\bigl(gl(n,\C)\bigr)$, its
representation in $\End(\C^n)$, and the
related RE algebra $\La_R$.
\section{The $GL(n,\C)$-case.}
\label{sGL}
From now on, we concentrate on the case $\g=gl(n,\C)$,
$G=GL(n,\C)$ and its representation in $\M=\End(\C^n)$.
Let us fix the Cartan subalgebra $\h\subset \g$ as the subspace of
diagonal matrices in $\M$.
As above, we denote by $\Pi_\g$, $\Pi_\g^\pm$ the sets of all, positive, and negative
roots with respect to $\h$.
It is customary, in the $gl(n,\C)$-case under consideration, to take
the trace pairing for the invariant scalar product $(\:\cdot\:,\:\cdot\:)$ on
$\g$. Let $e_\al,e_{-\al}, \al \in \Pi_\g^+$
be the root vectors  normalized to $(e_\al,e_{-\al}) = 1$.
Put $h_i$ to be diagonal matrix idempotents $h_i=e^i_i$, $i=1,\ldots, n$;
they  form  an orthonormal basis in $\h$.
The standard classical r-matrix and the invariant symmetric 2-tensor $\omega$
for $gl(n,\C)$  are
\be{}
r&=&\sum_{\al \in \Pi_\g^+}(e_{-\al}\tp e_\al - e_\al\tp e_{-\al}),\\
\omega &=& \sum_{i=1}^n  h_i\tp h_i + \sum_{\al \in \Pi_\g^+}
(e_{-\al}\tp e_\al + e_\al\tp e_{-\al}).
\label{omega}
\ee{}
Quantization of these data yields the Yang-Baxter operator
\be{}
\label{rsln}
R = q\sum_{i=1,\ldots ,n} e^i_i \tp e^i_i +
\sum_{i,j=1,\ldots ,n \atop i\not =j} e^i_i \tp e^j_j +
(q - q^{-1}) \sum_{i,k=1,\ldots ,n \atop i < k} e^k_i \tp e^i_k,
\ee{}
where $q = e^{h}$.
This is the image of the universal R-matrix
$\Ru\in\U^{\tp2}_h\bigl(gl(n,\C)\bigr)$.
The corresponding matrix $S$, defined by (\ref{S=PR}), satisfies
the Hecke condition
\be{}
S^2-(q-q^{-1})S=1\tp 1.
\label{Hecke}
\ee{}
We start our quantization programm\'e with computing the relevant
Poisson structures.

\subsection{RE Poisson structures on adjoint orbits of $GL(n,\C)$.}
\label{ssREPS}

Let $M$ be a right $G$-space and $r_M$ the bivector field on $M$
obtained from the classical r-matrix $r\in \wedge^2\g$ by
the group action.
Recall, \cite{DGS}, that if $\A_h(M)$ is a $\U_h(\g)$-equivariant quantization
on a right $G$-space $M$, the corresponding Poisson structure has the
form\footnote{We use slightly different definition of the quantum
group than in \cite{DGS}; this results in the different sign before $r_M$.}
\be{}
p=r_M+f,
\ee{}
where $r_M$ is the bivector field on  $M$ generated by the r-matrix
via the group action; $f$ is a $G$-invariant bivector field on $M$ whose Schouten
bracket is equal to
$$
\[f,f\] = -\[r_M,r_M\]=-\phi_M.
$$
Here $\phi \in \wedge^3 \g$ is an ad-invariant element,
see (\ref{phi}).
The group $G$ is equipped with the Poisson-Lie structure related to the classical
r-matrix $r$, \cite{Dr1}. Admissible Poisson brackets on $M$ are such that the
action $M\times G\to M$, where $M\times G$ is the Cartesian product of
Poisson manifolds, is a Poisson map.
The bivector $f$ defines a skew-symmetric bilinear operation on $\A(M)$ called
a $\phi$-bracket.
Specifically for the case $M=\End(V)$, where
the latter is considered as the right adjoint $G$-module,
the invariant part of the RE Poisson bracket is given
by the expression within parentheses in (\ref{rebr}).

Classification of $\phi$-brackets on semisimple orbits of semisimple Lie groups
was done in \cite{DGS} and \cite{Kar}. It this subsection, we compute those
obtained by restriction of the RE bracket (\ref{rebr}) to semisimple orbits
of the group $GL(n,\C)$ in $\End(\C^n)$.
Semisimple are orbits consisting of diagonalizable matrices.
They are characterized by eigenvalues and their multiplicities.
As abstract homogeneous spaces, they are specified by ordered sets of multiplicites which
fix the stabilizer subgroups $H\subset G$. Eigenvalues specify a Poisson structure
on the right coset space $H\backslash G$ by embedding it into the Poisson manifold $\End(\C^n)$.
We shall show that the RE Poisson structures on $H\backslash G$ form a variety
whose dimension coincides with the rank of the orbit.

\begin{lemma}
Let $X,Y$ be linear functions from $\M^*$ identified with $\M$ by the trace pairing.
The $\U(\g)$-invariant part of the reflection equation Poisson
bracket (\ref{rebr}) on $\M^\Omega$ is given by
\be{}
f(X,Y)(A) = (A^2,[X,Y])=\Tr(A^2[X,Y]), \quad A\in \M^{\Omega}.
\label{opb}
\ee{}
\end{lemma}
\begin{proof}
We identify $gl(n,\C)$ with $\M$, as well as the tangent space at the point
$A\in \M$. Right- and left-invariant
vector fields on $\M$ take the form $\xi^l \tr X=\xi X$, $X\tl\xi^r = X \xi $,
 $\xi\in \g$, $X \in \M^*\sim \M$.
Calculating the invariant part of the RE bracket (\ref{rebr})
with the invariant element $\omega$ from (\ref{omega}), we find
$$
\bigl(
(X \tl\omega^r_1)(\omega^l_2 \tr Y)-(\omega^l_1\tr X)(Y\tl\omega^r_2 )\bigr)|_A
=
(X\Omega_1,A)(\Omega_2Y,A)-(\Omega_1X,A)(Y\Omega_2,A).
$$
Here, we used the conventional notation with implicit summation, $\omega=\omega_1\tp \omega_2$
and $\Omega=\Omega_1\tp \Omega_2$.
The $\U\bigl(gl(n,\C)\bigr)$-invariant element $\Omega$ coincides with
the matrix permutation $P$, see (\ref{S=PR}).
Using the identity $(\Omega,X\tp Y)=(X,Y)$, which is valid for any $X,Y\in \M$,
we come to formula (\ref{opb}).
\end{proof}

Introduce notation $G_m=GL(m)$, $m\in \N$, and put
$G_{[n_1,\ldots,n_k]}=G_{n_1}\times\ldots \times G_{n_k}$, a Levi subgroup
in $G_n$. To every set $\{n_i\}$ of $k$ positive integers such that
$\sum_{i=1}^k n_i= n$ corresponds the right coset space
$O_{[n_1,\ldots,n_k]}=G_{[n_1,\ldots,n_k]}\backslash G_n$.
This becomes a  one-to-one correspondence between classes of
isomorphic homogeneous manifolds and sets $\{n_i\}$ provided they are
ordered.

An abstract homogeneous space, $O_{[n_1,\ldots,n_k]}$, is realized by orbits
in $\End(\C^n)$. Every such realization induces a Poisson structure
on it via restriction of the RE bracket (\ref{rebr}),
and different orbits give different Poisson structures.
Consider the direct sum decomposition $\C^n=\C^{n_1}\oplus\ldots \oplus\C^{n_k}$
and set $P_{n_i}\colon\C^n\to \C^{n_i}$ be the diagonal projector of rank $n_i$,
$i=1,\ldots,k$.
We define the orbit $O_{[n_1,\ldots,n_k;\la_1,\ldots,\la_k]}$
as that passing through the point $\sum_{i=1}^k \la_i P_{n_i}$, where
$\la_i$ are pairwise distinct
complex numbers. This correspondence between orbits and diagonal matrices
becomes one-to-one if we require some linear ordering among
those parameters $\la_i$ that correspond
to equal $n_i$. We choose the lexicographic odering on $\C$:
 $\la_1 \succ \la_2 $ $\iff$ $\Re\: \la_1> \Re \:\la_2$
or $\Re \:\la_1 = \Re\: \la_2$ and $\Im\: \la_1 >\:\Im \la_2$.
To summarize, the abstract homogeneous manifolds with
the RE Poisson structures  are in one-to-one correspondence
with ordered sets
of pairs $(n_i,\la_i)\in \N\times\C$
such that $n_i>0$ and $\sum_i n_i=n$. The ordering on $\N\times\C$
is defined as
\be{}
\label{order}
(n_1,\la_1) \succ (n_2,\la_2) \quad \Longleftrightarrow
\quad n_1> n_2 \quad \mbox{or}\quad  n_1 = n_2 \>\>\mbox{and}\>\>\la_1\succ\la_2.
\ee{}

Our further goal is to compute RE Poisson brackets, which
are distinguished by  their invariant parts,
on the abstract homogeneous manifolds  $O_{[n_1,\ldots,n_k]}$.
Let $\l\subset \g =gl(n,\C)$ be the Lie algebra of the sabilizer
$G_{[n_1,\ldots,n_k]}$ and  $\Pi_\l \subset\Pi_\g$ its
root system. The set of quasi-roots $\Pi_{\g/\l}$ corresponding
to a given Levi subalgebra $\l$ consists of  equivalence classes,
$\Pi_{\g/\l}=(\Pi_\g - \Pi_\l) \mod \Z\Pi_\l$, \cite{DGS}.
For any pair of quasi-roots $\bar\al$ and $\bar\bt$, it is possible to choose such representatives
$\al$ and $\bt$ that $\overline{\al+\bt} = \bar\al+\bar\bt$.
The root vectors
$e_\al$, $\al\in \bar\al \in \Pi_{\g/\l}$ form a basis
of the tangent space $\g/\l$ at the origin of $O_{[n_1,\ldots,n_k]}$.
When $\g=gl(n,\C)$, we can consider $\h^*$ as a subspace in
$\M$, using the invariant trace pairing. All roots of $\g$
are parameterized by pairs of integers $i,j=1,\dots,n$, $i\not=j$, and
can be represented by the elements $\al_{ij}=e^i_i-e^j_j$.

By definition, the orbit
$O_{[n_1,\ldots,n_k;\la_1\ldots,\la_k]}$ passes through
the point $\sum_{i=1}^k \la_i P_{n_i}$. The corresponding
quasi-roots are labeled by the pairs of integers
$i,j=1,\ldots,k-1$, $i\not =j$.
Every quasi-root $\bar \al_{ij}$ is the set of
elements $h_i-h_j$, where $h_i$ are diagonal matrix idempotents
of rank one such that $h_i P_{n_j} = \delta_{ij }h_i = P_{n_j} h_i$.

\begin{propn}
\label{REPS}
Consider an ordered set of $k$ positive integers $n_1\geq n_2\geq\ldots\geq n_k$.
The invariant part of RE Poisson bracket (\ref{rebr}) on the homogeneous space $O_{[n_1,\ldots,n_k]}$
is induced by the right-invariant bivector field
\be{}
\label{ip}
\sum_{\bar \al_{ij}\in \Pi^+_{\g/\l}}\frac{\la_i+\la_j}{\la_i-\la_j}\sum_{\al\in \bar \al_{ij}}
(e_\al\tp e_{-\al} - e_{-\al}\tp e_\al),
\ee{}
where $k$ complex numbers $\{\la_i\}$ form a decreasing sequence of
pairs $(n_i,\la_i)$, $i=1,\ldots,k$.
\end{propn}
\begin{proof}
As was shown in \cite{DGS}, an invariant bracket on
$O_{[n_1,\ldots,n_k]}$ must have the form
\be{}
\sum_{\bar \al_{ij}\in \Pi^+_{\g/\l}}c(\bar \al_{ij})\sum_{\al\in \bar \al_{ij}}
(e_\al\tp e_{-\al} - e_{-\al}\tp e_\al).
\label{ip1}
\ee{}
To find the coefficients $c(\bar \al_{ij})$, let us evaluate (\ref{ip1})
on the linear functions $e_{-\al}, e_{+\al} \in \M^*\simeq \M$, $\al\in\bar \al_{ij}$, at a diagonal matrix $A$:
$$
c(\bar \al_{ij})  e_{-\al}([A,e_{\al}]) e_{\al}([A,e_{-\al}])
=
c(\bar \al_{ij}) ([e_{\al},e_{-\al}],A) ([e_{-\al},e_{\al}],A)
=
-c(\bar \al_{ij}) \bigl({\al}(A)\bigr)^2.
$$
On the other hand,  formula (\ref{opb}) gives
$$
f(e_{-\al},e_{\al})|_A=(A^2,[e_{-\al},e_{\al}])=-\al(A^2).$$
Substituting $A=\sum_{i=1}^k \la_i P_{n_i}$, we obtain the coefficients $c(\bar \al_{ij})$:
$$
c(\bar \al_{ij})=\frac{\al(A^2)}{\bigl(\al(A)\bigr)^2}=\frac{\la^2_i-\la^2_j}{(\la_i-\la_j)^2}.
$$
\end{proof}
\noindent
Remark that the RE Poisson structures on $O_{[n_1,\ldots,n_k]}$
form a $k-1$-dimensional variety
because bivector (\ref{ip}) is stable under the dilation transformation
$\la_i\to \nu\la_i$, $i=1,\ldots,k$, $\nu\not= 0$.

\subsection{Characters of the RE algebra relative to the standard quantum group $\U_h\bigl(gl(n,\C)\bigr)$.}
\label{ssREAC}
In this subsection, we formulate the classification theorem for characters of
the algebra $\La_R$ associated with the representation
of the standard $\U_h\bigl(gl(n,\C)\bigr)$ in $\End(\C^n)$.
To do so, we need the following data.
\begin{definition}
\label{ap}
{\em An admissible pair} $(Y,\si)$ consists of a subset
$Y\subset I=\{1,\ldots ,n\}$ and a decreasing injective map
$\si\colon Y \to I$.
\end{definition}
\noindent
Clearly such a map is uniquely
determined by its image $\si(Y)$.
Introduce the subsets $Y_+=\{i\in Y| i>\si(i)\}$
and $Y_-=\{i\in Y| i<\si(i)\}$. Let
$b_-=\max\{Y_-\cup \si (Y_+)\}$ and  $b_+=\min\{Y_+\cup \si (Y_-)\}$.
Because $\si$ is a  decreasing map, one has
$b_-<b_+$.
\begin{thm}[\cite{M}]
\label{clth}
For the standard $gl(n,\C)$ R-matrix (\ref{rsln}), the
numerical solutions to the RE equation (\ref{nre}) fall into the following two classes.
\begin{enumerate}
\item Let $(Y,\si)$ be the admissible pair such that
$Y=[1,m]\cup[l+1,l+m]$, where $m$ and $l$ are non-negative integers
such that  $l+m\leq n$ and $m \leq l$; then $\si(i)=l+m+1-i$ for $i\in Y$.
Solutions of type A are gauge equivalent to
$$
A(l,m;\la,\mu)= \mu\sum_{i=1}^{m}e^i_i +
   \la\sum_{i=1}^{l}e^i_i +
   \sqrt{\la\mu}\sum_{i=1}^{m}(e^i_{\si(i)}-e_i^{\si(i)}),
$$
where $\mu,\la$ are arbitrary complex numbers.
The matrix $A(l,m;\la,\mu)$ has eigenvalues $\mu,\la,0$ with multiplicities
$m$, $l$, and $n-m-l$, correspondingly. It is semisimple if and only if
$\la\not =\mu$.
\item  Let $(Y,\si)$ be an admissible pair such that $\card(Y)\leq\frac{n}{2}$
and
$\si(Y) \cap Y=\emptyset$. Let $l$ be an integer
from the semiclosed interval $[b_-,b_+)$ and $\la \in \C$.
Solutions to the numerical RE of type B  are gauge equivalent
to
$$
B(Y,\si,l;\la)= \la\sum_{i=1}^{l}e^i_i +
   \sum_{i\in Y}e^i_{\si(i)}.
$$
The matrix $B(Y,\si,l;\la)$ has eigenvalue $\la$ and $0$ of
multiplicities $l$ and $n-l$. It is semisimple if and only if $\la\not =0$,
otherwise it is nilpotent of nilpotence degree  two.
\end{enumerate}
\end{thm}
\noindent
The layout of solutions to the numerical RE is as follows.
Generic RE matrix of type A is obtained  by
embedding $A(l,m;\la,\mu)|_{l+m=n}$ as the left top block extended with zeros to the entire matrix.
The matrix $A(l,m;\la,\mu)|_{l+m=n}$ itself has the form

\vspace{12pt}
\begin{picture}(200,100)
\put(70,45){$A(l,m;\la,\mu)|_{l+m=n}=$}
\put(300,45){.}

\put(183,-5){\line(0,1){110}}
\put(185,-5){\line(0,1){110}}
\put(295,-5){\line(0,1){110}}
\put(297,-5){\line(0,1){110}}
\put(224,34){\line(0,1){32}}
\put(256,34){\line(0,1){32}}
\put(224,34){\line(1,0){32}}
\put(224,66){\line(1,0){32}}

\put(195,91){$*$}
\put(188,81){$\mu+\la$}
\put(215,70){$*$}

\put(261,70){$*$}
\put(259,81){$\sqrt{\la\mu}$}
\put(282,91){$*$}

\put(194,3){$*$}
\put(187,14){$-\sqrt{\la\mu}$}
\put(215,24){$*$}

\put(230,54){$*$}
\put(237,47){$\la$}
\put(246,38){$*$}

\end{picture}
\vspace{12pt}

\noindent
The $(l-m)\times (l-m)$ square in the middle is located in the center of the matrix
and disappears when $l=m$.

Solutions of type B are decomposed into the direct sum
of matrices  $\la e^i_i$, $i\not \in Y\cup\si(Y), i<l$,
$\la e^i_i+ e^i_{\si(i)}$, ${i\in Y_-}$, and
$\la e^{\si(i)}_{\si(i)} +e^i_{\si(i)}$, ${i\in Y_+}$.

\subsection{Quantization of symmetric and bisymmetric orbits.}
\label{ssQSBO}
In this subsection, we apply the method of characters to quantization
of the RE bracket (\ref{rebr}) restricted to adjoint orbits
of $GL(n,\C)$ in $\End(\C^n)$.

\begin{propn}
The sub-variety $O_{[n_1,\ldots,n_k;\la_1,\ldots,\la_k]}\subset\M$
is determined by the system of equations
\be{}
(A-\la_1)\ldots (A-\la_k)&=&0,
\label{pc}\\
\Tr(A^m)&=& \sum_{i=1}^k n_i \;\la^m_i, \quad m=0,\ldots,k-1.
\label{trc}
\ee{}
\end{propn}
\begin{proof}
Reduction to the canonical Jordan form.
\end{proof}
Relations (\ref{pc}) and (\ref{trc}) may be generalized to
the quantum case, and one can try to use them for building
quantization of semisimple orbits.
However, the problem
is how to ensure the quotient by those relations be
a flat module over $\C[[h]]$.
The method of quantum character enables us to do that for certain
types of orbits.

Relation (\ref{pc}) makes sense in the algebra $\La_R$, because the
entries of any polynomial in generating matrix $L$ form a
$\U_h\bigl(gl(n,\C)\bigr)$-module
of the same type as
$L$ itself. There is a q-analog of the matrix trace, too.
Let $D$ be a matrix $D\in \M[[h]]$ such that the linear functional
$A\to\Tr(DA)$ on $\M[[h]]$ is invariant with respect to the action
$A\to \rho\bigl(\gm(x_{(1)})\bigr)A\rho(x_{(2)})$ of the quantum group
$\U_h\bigl(gl(n,\C)\bigr)$.
It is unique up to a scalar factor and we take it in the form
\be{}
\label{D}
D = \sum_{i=1}^{n} q^{-2i+2}e^i_i.
\ee{}
\begin{definition}
{\em Quantum trace} of a matrix $A$ with entries in
an associative algebra $\A$ is the element
\be{}
\label{PrTr}
\Tr_q(A) &=& \sum_{i=1}^{n} q^{-2i+2} A^i_i =\Tr(DA) \in \A.
\ee{}
\end{definition}
\noindent
When $\A=\La_R$ and $A=L$, the matrix of generators, this is an
invariant element belonging to the center of $\La_R$.
Moreover,  the  quantum  traces  $\Tr_q(L^k)$  of  the  powers in $L$ are
invariant and central as well,  \cite {AFS,KS}.
\begin{lemma}
\label{TrP}
Let $L$ be an RE matrix with coefficients
in an associative algebra $\A$. Suppose the matrix coefficients
$T^m_n\in \U^*_h(\g)$ commute with all $L^i_j$.
Then, the quantum trace is invariant under  the  similarity transformation
with the matrix $T$,
\be{}
\Tr_q  \bigl(\gm(T) L T\bigr)  & = &\Tr_q(L).
\label{Tr}
\ee{}
For any polynomial  ${\mathcal P}$ in one variable,
\be{}
{\mathcal P}\bigl(\gm(T) L T\bigr) & = & \gamma(T){\mathcal P}(L)T.
\label{P}
\ee{}
\end{lemma}
\begin{proof}
Formula (\ref{P}) is an immediate corollary of the equality
$\gm(T)= T^{-1}$. Verification of
(\ref{Tr}) is less simple and uses relations
(\ref{frt}) and (\ref{re}) in the algebras $\U^*_h(\g)$
and $\La_R$. The prove can be found, e.g., in
\cite{KS}.
\end{proof}
\begin{definition}[Quantum integers]
By {\em quantum integer} $\hat k$, $k\in \Z$, we mean the quantity
$\hat k =\frac{1-q^{-2n}}{1-q^{-2}}$.
\end{definition}
\noindent
Obviously, $\hat k$
is equal to the quantum trace,
$\hat k =\sum_{i=1}^{k} q^{-2i+2}$, of the unit endomorphism of the space $\C^k$,
for $k\in \N$.

Recall that we consider the set $\N \times \C$ lexicographically ordered,
see  (\ref{order}).
\begin{thm}
\label{sym}
Consider two pairs $(l,\la)$ and $(m,\mu)$ from
$\N \times \C$ and assume
$(l,\la) \succ (m,\mu) $, $l+m=n$.
The quotient of the algebra $\La_R$ by the relations
\be{}
(L-\la)(L-\mu)&=&0,
\label{QP}
\\
\Tr_q(L)&=&
\la\widehat {l} + \mu\widehat {m},
\label{Tr1}
\ee{}
is the $\U_h\bigl(gl(n,\C)\bigr)$-equivariant quantization of the
manifold $M=O_{[l,m]}$ with the Poisson bracket
$$
r_M+
\zeta\sum_{\al\in \bar \al_{12}}
(e_\al\tp e_{-\al} - e_{-\al}\tp e_\al),
$$
where $\zeta=\frac{\la+\mu}{\la-\mu}$.
\end{thm}
\begin{proof}
The Poisson structure on $O_{[l,m]}$ is induced by embedding  $O_{[l,m]}$
 in
$\M$ as the orbit  $O_{[l,m;\la,\mu]}$. The RE matrix
$A_h=A(l,m;\la,\mu)|_{l+m=n}$ does not depend on the deformation
parameter and belongs to $O_{[l,m;\la,\mu]}$. Applying Theorem \ref{thmquan}, we quantize this orbit by the
quantum character corresponding to $A_h$.
As a subalgebra in $\U^*_h\bigl(gl(n,\C)\bigr)$ the  algebra $\A_h(O_{[l,m]})$ is generated by
entries of the RE matrix $\gm(T)A T$, which fulfills
conditions (\ref{QP}) and (\ref{Tr1}), by  Lemma \ref{TrP}.
\end{proof}
\noindent
Note that if one of the eigenvalues $\mu, \la$ turns to zero, there are
several numerical RE matrices giving quantization of the same orbit.
We can take, e.g., the matrix $B(Y,\si,l;\la)$ for $A_h$, with
arbitrary set $Y$ such  that $\max\{Y_-\cup \si(Y_+)\}<l$. It has
eigenvalues $\la$ and $0$ of multiplicities $l$ and $n-l$.
For example, one can pick $Y=\emptyset$ and consider the RE matrix
$\la\sum_{i=1}^l e^i_j$. This solution to the numerical
RE and the corresponding quantization was built
in \cite{DM1}.
\begin{thm}
\label{bisym}
Consider two pairs $(l,\la)$ and $(m,\mu)$ from  $\N \times \C$
 and assume
$(l,\la) \succ (m,\mu)$, $n-(l+m)=k>0$.
The quotient of the algebra $\La_R$ by the relations
\be{}
L(L-\la)(L-\mu)&=&0,
\label{CP}
\\
\Tr_q(L)&=&
\la\widehat {l} + \mu\widehat {m},
\label{Tr11}
\\
\Tr_q(L^2)&=&
(\la+\mu)(\la\widehat {l} + \mu\widehat {m}) -\la\mu\: \widehat {l+m}
\label{Tr2}
\ee{}
is the $\U_h\bigl(gl(n,\C)\bigr)$-equivariant quantization
of the following manifolds:
\begin{enumerate}
\item
$M=O_{[l,m,k]}$ with the Poisson bracket
$$
r_M+
\zeta\sum_{\al\in \:\bar \al_{12}}
(e_\al\tp e_{-\al} - e_{-\al}\tp e_\al)+
\sum_{\al\in \:\bar \al_{23}\cup \:\bar \al_{13}}
(e_\al\tp e_{-\al} - e_{-\al}\tp e_\al),
$$
if $(l,\la) \succ (m,\mu)\succ (k,0)$,
\item
$M=O_{[l,k,m]}$ with the Poisson bracket
$$
r_M+
\zeta\sum_{\al\in\: \bar \al_{13}}
(e_\al\tp e_{-\al} - e_{-\al}\tp e_\al)+
\sum_{\al\in \:\bar \al_{12}\cup \:\bar \al_{32}}
(e_\al\tp e_{-\al} - e_{-\al}\tp e_\al),
$$
if $(l,\la)\succ (k,0)\succ (m,\mu)$,
\item
$M=O_{[k,l,m]}$ with the Poisson bracket
$$
r_M+
\zeta\sum_{\al\in\: \bar \al_{23}}
(e_\al\tp e_{-\al} - e_{-\al}\tp e_\al)+
\sum_{\al\in \:\bar \al_{21}\cup \:\bar \al_{31}}
(e_\al\tp e_{-\al} - e_{-\al}\tp e_\al),
$$
if $ (k,0)\succ (l,\la) \succ (m,\mu)$,
\end{enumerate}
where $\zeta=\frac{\la+\mu}{\la-\mu}$.
\end{thm}
\begin{proof}
As Poisson manifolds, all the three possibilities are
realized by the bisymmetric orbit passing through the point $A_h=A(l,m;\la,\mu)$ (cf. Theorem \ref{clth}),
which satisfies the RE equation and defines a character of
the RE algebra. By Theorem \ref{sym}, the quantization of the orbit with this character
is the quotient of the RE algebra and the subalgebra in
$\U^*_h\bigl(gl(n,\C)\bigr)$ generated by entries of  the RE matrix $\gamma(T)A_h T$.
Since the matrix $A_h $ satisfies conditions
(\ref{CP}--\ref{Tr2}), so does the matrix  $\gamma(T)A_h T$, by
Lemma \ref{TrP}.
\end{proof}
\begin{remark}
There is a one-parameter family of RE Poisson structures on
symmetric orbits, as follows from Proposition \ref{REPS},
 and their quantization is described by
Theorem \ref{sym}.  On bisymmetric orbits,
the RE Poisson brackets form a two-parameter family.
Theorem \ref{bisym} provides quantization for only
special one-parameter sub-families. This is a limitation
of the method of  characters, which gives those and only those
quantizations that can be  represented as subalgebras in
$\U^*_h(\g)$ and quotients of $\La_R$.
\end{remark}
\begin{remark}
Theorem \ref{clth} provides two classes of
non-semisimple numerical RE matrices obtained by the limits
$A(l,m;\la,\mu)|_{\la\to\mu}$ and  $B(Y,\si,l;\la)|_{\la\to0}$.
They belong to orbits
that are limits of the semisimple ones.
The orbit passing through $B(Y,\si,l;\la)|_{\la=0}$
is nilpotent. Using Theorem \ref{thmquan}, one can quantize
all nilpotent orbits of nilpotence degree two.
\end{remark}

\subsection{Quantizing fiber bundle $O_{[l,m,k]}\to O_{[l+m,k]}$.}
\label{ssQFB}

In this subsection we consider  the following problem.
The group embedding
$G_{[l,m,k]}  \to G_{[l+m,k]}$
defines a map of the right coset
spaces,
\be{}
O_{[l,m,k]} \to O_{[l+m,k]},
\label{fiber}
\ee{}
which is a bundle
with the fiber
$O_{[l,m]}=  G_{[l,m,k]}\backslash G_{[l+m,k]}\simeq G_{[l,m]}\backslash G_{l+m}$.
Map (\ref{fiber}) is equivariant with respect to the action of the
group $G_{l+m+k}$ and the bundle is homogeneous.
Suppose the total space and the base of the bundle are Poisson manifolds
and map (\ref{fiber}) is Poisson.
The problem is to quantize the diagram
\be{}
\A(O_{[l,m,k]}) \leftarrow \A(O_{[l+m,k]}),
\label{fiber0}
\ee{}
i.e., to build $\U_h\bigl(gl(l+m+k)\bigr)$-equivariant quantizations  of the function algebras and
an equivariant monomorphism
\be{}
 \A_h(O_{[l,m,k]}) \leftarrow \A_h(O_{[l+m,k]}).
\label{qfiber}
\ee{}
The spaces $O_{[l,m,k]}$  and $O_{[l+m,k]}$ can be realized as
the  orbits  $O_{[l,m;k;\mu,\la,0]}$ and  $O_{[l+m,k;-\la\mu,0]}$
(here we do not assume ordering (\ref{order}))
with the induced RE Poisson structures. We will show
that their quantizations built in the previous subsection
admit an equivariant morphism, a quantization
of (\ref{fiber0}). This will imply that the projection
$O_{[l,m;k;\la,\mu,0]}\to O_{[l+m,k;-\la\mu,0]}$ is a Poisson map,
because of the flatness of the quantizations.

\begin{remark}
The fiber over the origin $A(l,m;\la,\mu)\in O_{[l,m,k;\la,\mu,0]}$
is realized as a Poisson
manifold $O_{[l,m;\la,\mu]}$. It can be shown that its embedding to
$O_{[l,m,k;\la,\mu,0]}$ is a Poisson map as well and this map
can be lifted to a homomorphism of quantized algebras
that is equivariant with respect to the quantum group embedding
$\U_h\bigl(gl(l+m)\bigr)\to \U_h\bigr(gl(l+m+k)\bigr)$.
The proof of this statement is rather straightforward, and we do not
concentrate on this subject here.
\end{remark}

Before formulating the main result of this subsection, let us prove an algebraic
statement.
Let ${\mathcal P}(x)$ be a polynomial
in one variable with coefficients in a commutative ring $\K$. Suppose
${\mathcal P}(0)=0$. Fix two scalars  $\al,\bt \in \K$ and consider an
associative unital algebra
$\A(e,s)$ over
$\K$ generated by the  elements $\{e,s\}$ subject to relations
\be{}
eses=sese \quad\mbox{("reflection equation")},\quad  s^2-\bt s=1
 \quad\mbox{("Hecke condition")}.
\label{REHc}
\ee{}
\begin{lemma}
\label{alg}
The correspondence $s\to s, e\to {\mathcal P}(e)$
gives a homomorphism of the algebra $\A(e,s)$
to the quotient algebra $\A(e,s)/(e{\mathcal P}(e)-\al e)$.
This homomorphism is factored through the ideal $(e^2-\al e)$.
\end{lemma}
\begin{proof}
First let us note that the last statement is an immediate corollary
of the condition ${\mathcal P}(0)=0$.
We should verify that, modulo the ideal $(e{\mathcal P}(e)-\al e)$,
the following relation holds true in the algebra $\A(e,s)$:
\be{}
\label{weaker}
{\mathcal P}(e)s{\mathcal P}(e)s = s{\mathcal P}(e)s {\mathcal P}(e).
\ee{}
We will prove a stronger assertion; namely,
for any $m=1,2,\ldots $, the identity
\be{}
\label{stronger}
{\mathcal P}(e)se^ms = se^ms {\mathcal P}(e)
\ee{}
is valid modulo the ideal $(e{\mathcal P}(e)-\al e)$. This will imply
(\ref{weaker}), because ${\mathcal P}(0)=0$ by the hypothesis of the lemma.
We assume $\bt \not =0$ since otherwise
$se^ms$ can be replaced by $(ses)^m$, and the prove becomes immediate.
For $m=1$ this is a consequence of the "reflection equation" relation
(\ref{REHc}).
Suppose (\ref{stronger}) is proven for some integer $m=l\geq 1$.
Using the "Hecke condition" (\ref{REHc}),
we rewrite (\ref{stronger}) for $m=l+1$ as
$${\mathcal P}(e)se^l(s^2-\bt s)es = se^l(s^2-\bt s)es {\mathcal P}(e).$$
The terms without $\bt$ compensate each other, by the induction assumption.
The problem reduces to checking the equality
$${\mathcal P}(e)se^l ses = se^l s es {\mathcal P}(e).$$
Employing the induction assumption, rewrite
the left-hand side as
$${\mathcal P}(e)se^l ses = se^l s{\mathcal P}(e)es = \al se^l s es.$$
The right hand side can be transformed as
$$
se^l s es {\mathcal P}(e) = se^{l-1} e s es {\mathcal P}(e)=
se^{l-1}  s ese {\mathcal P}(e)=\al se^{l-1}  s ese = \al se^{l-1}  e s es=\al se^l s es.
$$
Here, we used the "reflection equation"  (\ref{REHc}).
\end{proof}
\begin{thm}
Let $E$ and $B$ be the RE matrices whose entries generate the algebras
$\A_h\bigl(O_{[l,m,k;\la,\mu,0]}\bigr)$ and $\A_h\bigl(O_{[l+m,k;-\la\mu,0]}\bigr)$,
respectively.
The matrix correspondence $\pi(B) = E^2 - (\la+\mu)E$ is extended to the
$\U_h\bigl(gl(l+m+k)\bigr)$-equivariant algebra morphism
\be{}
\A_h\bigl(O_{[l,m,k;\la,\mu,0]}\bigr)\leftarrow  \A_h\bigl(O_{[l+m,k;-\la\mu,0]}\bigr).
\label{qofiber}
\ee{}
\end{thm}
\begin{proof}
The matrix $\pi(B)$ satisfies the equations $\pi(B)(\pi(B)+\la\mu)=0$ and
\be{}
\Tr_q\bigl(\pi(B)\bigr)&=&
(\la+\mu)(\la\widehat {l} + \mu\widehat {m}) -\la\mu \:\widehat {l+m}
-(\la+\mu)(\la\widehat {l} + \mu\widehat {m})
\n\\&=&
-\la \mu\:\widehat {l+m}
\n.
\ee{}
Taking into account Theorem \ref{sym}, it remains to show that $\pi(B)$
is an RE matrix.
The matrix $E$ fulfills polynomial relation (\ref{CP}), and
the braid matrix $S$ matrix satisfies the Hecke  condition
(\ref{Hecke}). It remains to apply
Lemma \ref{alg}, setting $e=E_2$, $s=S$,
${\mathcal P}=\pi$, $\al=-\la\mu$, $\bt=q-q^{-1}$.
\end{proof}
\subsection{Quantizing the Kirillov-Kostant-Souriau bracket on symmetric orbits.}
\label{QKKS}
It is known that there is a two-parameter quantization $\La_{h,t}$ on
$\End(\C^n)$ that is equivariant with respect to the
adjoint action of $\U_h\bigl(gl(n,\C)\big)$, \cite{D1}. The corresponding
Poisson structure is obtained from (\ref{rebr}) by adding the
Poisson-Lie bracket with arbitrary overall factor $t$.
The algebra $\La_{h,t}$ can be obtained from $\La_R$ by
the substitution $L=E+\frac{t }{1-q^{-2}}$, $q=e^h$, of the matrix of generators.
Relations (\ref{re}) go over into
\be{}
\label{cuea}
SE_2SE_2-E_2SE_2S=qt\:[E_2,S].
\ee{}
This is true for any matrix $S$ satisfying the Hecke condition
(\ref{Hecke}).
In the limit $h\to 0$, the matrix $S$ tends to the permutation
operator $P$ on  $\C^n\tp \C^n$ and relations (\ref{cuea}) turns into
those of the classical universal enveloping algebra
$\U\bigl(gl(n,\C)\big)[t]$. Indeed, the substitution $P\to S$
in (\ref{cuea}) gives explicitly
\be{}
\label{clcr}
E^i_j E^m_n - E^m_n E^i_j = t (\delta^m_j E^i_n-\delta^i_n E^m_j).
\ee{}
The algebra $\U\bigl(gl(n,\C)\big)[t]$
is a quantization of
the Poisson-Lie bracket on $\End^*(\C^n)$ identified with
$gl(n,\C)$ by the invariant trace pairing.
Its restriction to orbits is the Kirillov-Kostant-Souriau  bracket.
\begin{thm}
\label{2param}
Let $\mu_1$,  $\mu_2$ be two distinct complex numbers
and $n_1$, $n_2$ two positive integers such that $n_1+n_2=n$.
Let $\{E^i_j\}_{i,j=1}^n$ be the set of generators of
the algebra $\La_{h,t}$ subject to relations
(\ref{cuea}).
The quotient of $\La_{h,t}$ by the
relations
\be{}
(E-\mu_1)(E-\mu_2)&=&0,
\label{q-tsym1}
\\
\Tr_q(E)&=&\hat{n}_1\mu_1+\hat{n}_2\mu_2 + t \hat{n}_1\hat{n}_2.
\label{q-tsym2}
\ee{}
is a two-parameter quantization
on the symmetric orbit
$O_{[n_1,n_2;\mu_1,\mu_2]}$. It is equivariant
with respect to the adjoint action of $\U_h\bigl(gl(n,\C)\bigr)$.
\end{thm}
\begin{proof}
Consider the quantization of $O_{[l,m;\la,\mu]}$ as
a quotient of the RE algebra $\La_R$, along the line of Theorem \ref{sym}.
The substitution
$$
L=E+\frac{t}{1-q^{-2}},\quad \la=\mu_1+\frac{t}{1-q^{-2}},\quad
\mu=\mu_2+\frac{t}{1-q^{-2}},\quad
l=n_1,\quad m=n_2
$$
transforms relations (\ref{re}) into (\ref{cuea})
while (\ref{QP}--\ref{Tr1}) into
(\ref{q-tsym1}--\ref{q-tsym2}).
\end{proof}
\begin{corollary}
\label{QKKSbr}
Let $\mu_1$,  $\mu_2$ be two distinct complex numbers
and $n_1$, $n_2$ two positive integers such that $n_1+n_2=n$.
Let $\{E^i_j\}_{i,j=1}^n$ be the set of generators of
the algebra $\U\bigl(gl(n,\C)\big)[t]$ subject to relations
(\ref{clcr}).
The quotient of $\U\bigl(gl(n,\C)\big)[t]$ with respect to  the
ideal generated by
\be{}
\label{KKSsym1}
(E-\mu_1)(E-\mu_2)&=&0,\\
\label{KKSsym2}
\Tr(E)&=&n_1\mu_1+n_2\mu_2  + t n_1n_2,
\ee{}
is the quantization of the KKS bracket
on the symmetric orbit
$O_{[n_1,n_2;\mu_1,\mu_2]}$. It is equivariant
with respect to the adjoint action of $\U\bigl(gl(n,\C)\bigr)$.
\end{corollary}
\begin{proof}
Taking the limit $h\to0$ in the two-parameter quantization
of Theorem \ref{2param}.
\end{proof}
\begin{remark}
Theorem  \ref{2param} gives a two-parameter
generalization  of
the quantum sphere, \cite{GS}, to symmetric orbits. The explicit description
of the quantized KKS bracket in Corollary \ref{QKKSbr}
is an especially interesting result that could not
be otherwise obtained than extending deformation
to the q-domain. This is a remarkable application
of the quantum group theory.
\end{remark}
\begin{example}[Complex sphere]
Now we illustrate  Corollary \ref{QKKSbr} on
$\U\bigl(gl(2,\C)\bigr)$-equivariant quantization
of the complex sphere $O_{[\mu_1,\mu_2;1,1]}\subset\End(2,\C)$.
Our goal is to demonstrate on this simple example that
the system of conditions (\ref{clcr}), (\ref{KKSsym1}), and (\ref{KKSsym2})
is self-consistent.
It is known that $O_{[\mu_1,\mu_2;1,1]}$ is
specified, as a maximal orbit, by values of two invariant functions, the
traces of a matrix and its square.
We shall show that matrix
equation  (\ref{KKSsym1}) boils down to a condition on
the second Casimir of $\U\bigl(gl(2,\C)\bigr)$ only
when the first Casimir is fixed as in (\ref{KKSsym2}).

Consider the $2\times2$-matrix $E=||E^i_j||$
of generators of  $\U\bigl((2,\C)\bigr)[t]$
obeying the commutation relations
\be{}
\begin{array}{llrcllr}
[E^1_1, E^1_2]&=&t E^1_2,& &
[E^1_1, E^2_1]&=&-t E^2_1,\\[6pt]
[E^2_2, E^2_1]&=&t E^2_1,& &
[E^2_2, E^1_2]&=&-t E^1_2,\\[6pt]
[E^1_1, E^2_2]&=&0,& &
[E^1_2, E^2_1]&=&t (E^1_1-E^2_2),
\end{array}{}
\label{gl2cr}
\ee{}
which are obtained by specialization of  system
(\ref{clcr}) to the $gl(2,\C)$-case.
Put $n_1=n_2=1$, $\si_1=\mu_1+\mu_2$, $\si_2=\mu_1\mu_2$
in (\ref{KKSsym1}--\ref{KKSsym2})
and write (\ref{KKSsym1}) out explicitly:
\be{}
(E^1_1)^2  +  E^2_1 E^1_2 - \si_1 E^1_1 + \si_2&=&0,
\label{e1}\\
(E^2_2)^2  +  E^1_2 E^2_1 - \si_1 E^2_2 + \si_2&=&0,
\label{e2}\\
E^1_2 E^1_1  +  E^2_2 E^1_2 - \si_1 E^1_2 &=&0,
\label{e3}\\
E^1_1 E^2_1  +  E^2_1 E^2_2 - \si_1 E^2_1 &=&0
\label{e4}.
\ee{}
Equations (\ref{e1}--\ref{e4}) are equivalent to
the system
\be{}
(E^1_1)^2 +(E^2_2)^2  +  E^1_2 E^2_1  +  E^2_2 E^2_1
- \si_1 (E^1_1+E^2_2) + 2\si_2&=&0,
\label{e1'}\\
(E^1_1)^2 - (E^2_2)^2 - t (E^1_1-E^2_2)   - \si_1 (E^1_1-E^2_2) &=&0
\label{e2'},\\
E^1_1 E^1_2  +  E^2_2 E^1_2 - (\si_1 +t) E^1_2 &=&0
\label{e3'},\\
E^1_1 E^2_1  +  E^2_2 E^2_1 - (\si_1 +t) E^2_1 &=&0
\label{e4'}.
\ee{}
To obtain it, we pulled the diagonal generators to the left in
(\ref{e3}--\ref{e4}) and took the sum and difference
of (\ref{e1}--\ref{e2}) using the commutation relations
(\ref{gl2cr}).

Equations (\ref{e2'}--\ref{e4'}) are satisfied modulo the condition
\be{}
E^1_1+E^2_2-(\si_1+t)&=&0
\label{e5'}.
\ee{}
which is the specialization of (\ref{KKSsym2}) for $O_{[\mu_1,\mu_2;1,1]}$.
Equations (\ref{e1'}) and (\ref{e5'}) can be rewritten
in terms of the elements $\Tr(E)=\sum_i E^i_i$,
$\Tr(E^2)=\sum_{i,j} E^i_j E^j_i$ that generate the center
of $\U\bigl(gl(2,\C)\bigr)$.
The resulting relations for the quantized  orbit
$O_{[\mu_1,\mu_2;1,1]}$
are
\be{}
\Tr(E^2)&=& \si_1 (\si_1+t) - 2\si_2,
\label{e1''}\\
\Tr(E)&=&(\si_1+t)
\label{e5''}.
\ee{}
Together with (\ref{gl2cr}), this is a
$\U\bigl((2,\C)\bigr)$-equivariant quantization of the KKS bracket
on the two-dimensional complex sphere.
\end{example}
\section{Conclusion.}
The method of quantum characters formulated in this paper is designed for building
$\U_h(\g)$-equivariant quantizations on a $G$-manifold
$M$ that are representable as subalgebras
in $\U^*_h(\g)$ and quotients of $\A_h(M)$. Despite
its simplicity, it allows to obtain new and interesting results,
for example, the two-parameter quantization
on semisimple orbits.
Analyzing the quantizations built within the present approach, one may
come to the following
conclusion. The two-parameter quantization on a semisimple
orbit of $GL(n,\C)$ may be sought for in the form of
a matrix polynomial equation on
the generators $E^i_j\in  \La_{h,t}$ with
additional conditions on the quantum traces $\Tr_q(E^k)$, $k\in \N$.
This conjecture turns out to be true.
The proof is based on a different technique than that used in the present
paper. It is the subject of our forthcoming publication, \cite{DM3}, as well
as the explicit equations defining quantized semisimple
orbits of $GL(n,\C)$, including the special case of the KKS bracket.

\vspace{0.3cm}
\noindent
{\large\bf Acknowledments.}
We thank Steven Shnider for numerous valuable discussions. We are grateful
to the referee for his remarks, which helped us to improve the manuscript.

\bigskip
e-mail: donin@macs.biu.ac.il\\
e-mail: mudrova@macs.biu.ac.il;

\end{document}